\documentstyle [12pt]{article}
\newfont{\bbb} {msbm10}
\newcommand{\Bbb}[1]{\mbox{\bbb#1}}
\newcommand{\R}{\Bbb{R}}

\newcommand{\N}{\Bbb{N}}

\newcommand{\Z}{\Bbb{Z}}

\newcommand{\ra}{\rightarrow}
\begin{document}
\title{Cocompact CAT(0) Spaces are Almost Geodesically Complete}
\author{Pedro Ontaneda}
\date{}
\maketitle

Let $M$ be a Hadamard manifold, that is, a complete simply connected 
riemannian manifold with
non-positive sectional curvatures. Then every geodesic segment
$\alpha :[0,a]\rightarrow M$ from $\alpha (0)$ to $\alpha (a)$ can be 
extended
to a geodesic ray $\alpha : [0,\infty )\rightarrow M$. 
We say then that the Hadamard manifold $M$ is geodesically complete.
Note that, in this case, all geodesic rays are proper maps.

CAT(0) spaces are generalizations of Hadamard manifolds.
For a CAT(0) space $X$, all geodesic rays $\alpha : [0,\infty )
\rightarrow X$ are proper maps but, in 
general, $X$ is not geodesically complete. The following definition of almost 
geodesic completeness 
was suggested by M. Mihalik:
\vspace{.1in}

{\it A geodesic space $X$, with metric d, is almost geodesically complete if there is a 
constant $C$
such that  for every $p,q\in X$ there is a geodesic ray $\alpha :
[0,\infty)\rightarrow X$, $\alpha (0)=p$, and $d(q,\alpha)\leq C$.}
\vspace{.1in}

In general CAT(0) and hyperbolic spaces are not almost geodesically 
complete.
For instance, the non-negative reals with a line segment of length $N$ attached
at the integer $N$ for all $N>0$ is not almost
geodesically complete. In the presence
of cocompact group actions the situation changes. 

Suppose $G $ is the Cayley graph of an infinite word hyperbolic group.
If $q$ is a vertex of $G$ then there is a geodesic line $l$ containing
$q$. For $p\in G$ let $r_1$ and $r_2$ be geodesic rays from $p$ to the two limit
points of $l$. This forms an ideal $\delta$-thin  triangle and so 
either
$r_1$ or $r_2$ must pass within $\delta$ of $q$. 
Hence $G$ is almost geodesically complete.

This basic fact about word hyperbolic groups is used extensively in the 
literature (see for example \cite{BM}) and Mihalik conjectured an 
analogous result should be true for CAT(0) groups, that is, groups acting 
cocompactly by isometries on CAT(0) spaces. For general CAT(0) spaces there 
are no thin triangles, hence the argument used for word hyperbolic 
groups above does not work. 
\vspace{.1in}

In this paper we prove that, under certain conditions, cocompact CAT(0)
spaces are almost geodesically complete. In the following statements
$H^{i}_{c}(X)$ denotes cohomology with integer coefficients and 
compact supports. Also, a metric space is  {\it proper} if all
closed balls are compact.\vspace{.1in}

Our main results are 
\vspace{.3in}

{\bf Theorem A.} {\it Let $X$ be a noncompact proper CAT(0) space on 
which $\Gamma$ acts cocompactly
by isometries. If $H^{i}_{c}(X)\neq 0$, for some $i$,
then $X$ is almost geodesically complete.}
\vspace{.3in}

{\bf Theorem B.} {\it Let $X$ be a noncompact proper CAT(0) space on 
which $\Gamma$ acts cocompactly by isometries with discrete orbits.
Then $X$ is almost geodesically complete.}
\vspace{.3in}

Theorem B follows from theorem A and the following two propositions.
\vspace{.2in}

{\bf Proposition A.} {\it Let $X$ be a proper CAT(0) space on which 
$\Gamma$ acts cocompactly by isometries  with discrete orbits.
 Then $X$ is properly $\Gamma$-homotopy equivalent
 to a $\Gamma$-finite $\Gamma$-simplicial complex $K$.}
\vspace{.2in}

{\bf Proposition B.} {\it Let $K$ be a locally finite
contractible simplicial complex which admits a cocompact simplicial action.
Then $H^{i}_{c}(K)\neq 0$, for some $i$.}
\vspace{.2in}

{\bf Proof of Theorem B from Theorem A and Propositions A and B.}
Let $X$ be a noncompact proper CAT(0) space on 
which $\Gamma$ acts cocompactly, by isometries with discrete orbits.
By Proposition A,  $X$ is properly $\Gamma$-homotopy equivalent
to a $\Gamma$-finite $\Gamma$-simplicial complex $K$.
Since every CAT(0) space is contractible, we have that $K$ is also
contractible. Hence, by Proposition A, $H^{i}_{c}(X)=H^{i}_{c}(K)\neq 0$.
We can now apply Theorem A and conclude that $X$ is almost geodesically complete.
\vspace{.2in}

It was suggested by R. Geoghegan that proposition A above could 
be used to
prove that that the boundary $\partial X$ of a $\Gamma$-cocompact 
CAT(0)
space $X$ is a shape invariant of the $\Gamma$ action. The next theorem 
shows that in fact this is true. \vspace{.3in}

{\bf Theorem C.} {\it Let $X$ and $Y$ be proper CAT(0) spaces on which 
$\Gamma$ acts cocompactly by isometries with discrete orbits. If 
$X$ and $Y$ are $\Gamma$-homotopy
equivalent then $\partial X$ and $\partial Y$ are shape equivalent.}
\vspace{.3in}

{\bf Corollary A.} {\it Let $X$ and $Y$ be proper CAT(0) spaces on which 
$\Gamma$ acts cocompactly by isometries with discrete orbits.
If the actions have the same isotropy, i.e. if 
$\{\, G<\Gamma\, :\, X^{G}\neq\emptyset\,\}=\{\,G<\Gamma\, :\, 
Y^{G}\neq\emptyset\,\}$
then $\partial X$ and $\partial Y$ are shape equivalent.}
\vspace{.2in}

We say that a group acts on a space with finite isotropy
if all isotropy groups are finite.
\vspace{.2in}

{\bf Corollary B.} {\it Let $X$ and $Y$ be proper CAT(0) spaces on which 
$\Gamma$ acts cocompactly by isometries, with discrete orbits and 
finite isotropy.
Then $\partial X$ and $\partial Y$ are shape equivalent.}
\vspace{.3in}

It is known that if we assume $\Gamma$ in theorem C to be hyperbolic, 
then in fact
$\partial X$ and $\partial Y$ are homeomorphic (see \cite{G}), but in 
general
$\partial X$ and $\partial Y$ do not have to be homeomorphic (see 
\cite{CKl}).

Here is a short outline of the paper. In section 1 we recall some 
definitions and a lemma. In section 2 we prove proposition A,
theorem C and its corollaries. 
In section 3 we prove theorem A. Finally, in section 4 we prove proposition B.

We are grateful to Mike Mihalik for suggesting the problem and for his 
helpful comments, and to Dan Farley for pointing out some improvements 
in the exposition.
The idea of using nerves in the proof of proposition A was suggested 
by Stratos Prassidis
and also by the work of Bieri and Geoghegan \cite{BG}. We are grateful 
to them. We thank the referee for its useful remarks.
\vspace{.3in}

{\bf 1. Definitions.}
\vspace{.2in}

For the definitions and basic facts about geodesics and CAT(0) spaces see, 
for instance, \cite{B} or \cite{BH}. Throughout this paper all CAT(0) spaces
are assumed to be {\it complete}  metric spaces.
Note that this assumption does not affect the statements of our main results since
we always take our CAT(0) spaces to be proper. Recall that a metric space
is {\it proper} if all balls are compact.

We say that a group $\Gamma$ {\it acts 
cocompactly} on a space $X$ if there is a compact subset $C$ of $X$ such that 
$X=\bigcup_{\gamma\in \Gamma}\gamma C$.

If $\Gamma$ acts on a space $X$ define the {\it isotropy groups}
$$\Gamma^{x}=\{\gamma\in\Gamma\,\,:\,\, \gamma x=x\,\}$$ and if 
$A\subset\Gamma$, define the subgroups 
$$\Gamma^{A}=\{\gamma\in\Gamma\,\,:\,\, \gamma \,\, fixes\,\, A\,\, 
pointwise\,\}$$
and $$\Gamma^{(A)}=\{\gamma\in\Gamma\,\,:\,\, \gamma A=A\,\}$$
Define also, 
for $G\subset\Gamma$, the {\it fixed point set}
 $$X^{G}=\{x\in X\,\, : \,\, g x=x,\,\, g\in G\}$$
Recall that if $X$ is a CAT(0) space and $\Gamma$ acts by isometries on 
$X$, then $X^{G}$ is convex, hence contractible.
Also,
note that $X^{G}\neq \phi$ iff $G\subset \Gamma^{x}$, for some
$x\in X$.

Let $\cal{C}$ be a collection of subgroups of $\Gamma$. We say that
$X$ is $\cal{C}$-{\it free} if $\Gamma^{x}\in \cal{C}$, for all $x\in X$. 
Also, we say that $X$ is $\cal{C}$-{\it contractible} if $X^{G}$ is non-empty
and contractible, for $G\in\cal{C}$. If $\Gamma$ acts cellularly
on a CW-complex $K$, we say that $K$ is a {\it universal} $(\Gamma , 
\cal{C})$-{\it complex}
 if $K$ is $\cal{C}$-free, $\cal{C}$-contractible and $\Gamma^{(\sigma 
)}=\Gamma^{\sigma}$, for all cells $\sigma$  of $K$.
\vspace{.2in}

{\bf Lemma 1.1.} {\it Every CAT(0) space is an AR.}
\vspace{.1in}

The proof follows from \cite{Hu}, IV.4.1. by taking a refinement of 
the covering $\alpha$ consisting of convex subsets (for instance balls).
Theorem IV.1.2 of \cite{Hu} also works.
\vspace{1.2in}

{\bf 2. Proof of Proposition A,  Theorem C and the Corollaries.}
\vspace{.2in}

{\bf 2.1. Proposition A.} {\it Let $X$ be a proper CAT(0) space on which 
$\Gamma$ acts cocompactly by isometries with discrete orbits.
 Then $X$ is $\Gamma$-homotopy equivalent to a $\Gamma$-finite 
$\Gamma$-simplicial complex $K$.}
\vspace{.2in}

{\bf Proof.}
Since the orbits $\Gamma x$  are discrete, we have that for every 
$x\in  X$ there is a closed ball $B_{x}(r_{x})$, with center $x$ and 
radius $r_{x}>0$, such that $B_{x}(r_{x})\cap(\Gamma x)=\{x\}$. Note that this 
implies that  for every $\gamma\in \Gamma$, either 
$B_{x}(\frac{r_{x}}{2})\cap\gamma B_{x}(\frac{r_{x}}{2})=\emptyset$ or $\gamma 
x= x$, and this 
last case implies
$B_{x}(\frac{r_{x}}{2})=\gamma B_{x}(\frac{r_{x}}{2})$.
Also, because $X$ is proper, these balls are compact.

Now, the action is cocompact, thus there is a finite  collection 
${\cal{ V}}$
of  balls $B_{x}(\frac{r_{x}}{4})$,  
 such that 
 $\bigcup_{\gamma\in\Gamma ,\,V\in{\cal{V}}}\gamma V=X$.
 Define ${\cal{U}}=\Gamma
{\cal{V}}=\{\gamma V\, :\,\,\gamma\in\Gamma,\,\, V\in {\cal{V}}\}$.

We show that every ball $U\in {\cal{U}}$ intersects only a finite 
number of elements in $\cal U$. 
Suppose not, then there are $V_{1}=B_{x}(\frac{r_{x}}{4})$, 
$V_{2}=B_{y}(\frac{r_{y}}{4})$ and a sequence $\{\gamma_{i}\}\subset\Gamma$, such 
that $\gamma_{i}V_{1}$ intersects
$V_{2}$, and the $\gamma_{i}V_{1}$'s are all different.
We may assume that $r_{x}\geq r_{y}$. Thus, if 
$\gamma_{i}B_{x}(\frac{r_{x}}{4})$ and $\gamma_{j}B_{x}(\frac{r_{x}}{4})$ 
both intersect $V_{2}$, then $\{ y\}\in \gamma_{i}B_{x}(\frac{r_{x}}{2})\cap
\gamma_{j}B_{x}(\frac{r_{x}}{2})$. Therefore
$\gamma_{i}B_{x}(\frac{r_{x}}{2})$ intersects  $\gamma_{j}B_{x}(\frac{r_{x}}{2})$,  
for all $i,j$. Hence $(\gamma_{j})^{-1}
\gamma_{i}B_{x}(\frac{r_{x}}{2})$ intersects $B_{x}(\frac{r_{x}}{2})$
and we mention before that this implies 
$(\gamma_{j})^{-1}\gamma_{i}x=x$. Follows that the 
$\gamma_{i}V_{1}$'s are all equal, a contradiction.
This proves that every ball $U\in\cal{U}$ intersects only a finite 
number of elements in $\cal{U}$.

Denote by $K$ the nerve of the covering $\cal U$. Recall that $K$ is 
the complex that has one vertex for each element of $\cal U$, and 
$\{U_{0},...,U_{k}\}$
 forms a simplex $\, <U_{0},...,U_{k}\, >\,$ if $U_{0}\cap ...\cap U_{k}$ is 
non-empty. Note that, since every element in $\cal{U}$ intersects only 
a finite number of elements in $\cal U$,  $K$ is locally finite. 
Remark that $\Gamma$ acts on $K$, simplicially by $\gamma 
<U_{0},...,U_{k}>\, =\, <\gamma U_{0},...,\gamma U_{k}>$. It follows that 
$K$, with this action, is $\Gamma$-finite and cocompact; hence $K$ is finite dimensional.
\vspace{.2in}

{\bf Claim 1.} {\it $\Gamma^{(\sigma )}=\Gamma^{\sigma}$, where 
$\sigma$ is a simplex of $K$.}
\vspace{.1in}

Let $\{U_{0},...,U_{k}\}$ be the vertices of $\sigma$ in $K$. For 
$\gamma\in \Gamma^{(\sigma )}=\{\gamma\in\Gamma :\,\, \gamma\sigma=\sigma\}$
we have that $\gamma (U_{0}\cap...\cap U_{k})
=U_{0}\cap...\cap U_{k}\neq\emptyset$.
 Hence $\gamma U_{i}$ intersects $U_{i}$, which means, as before, that 
$\gamma U_{i}=U_{i}$. Thus $\gamma$  
 fixes $\sigma$ pointwise and $\Gamma^{(\sigma )}=\Gamma^{\sigma}$. 
This proves claim 1.
\vspace{.2in}

{\bf Claim 2. } {\it $K^{G}$ is contractible, for  $G\subset\Gamma$, 
and $K^{G}\neq \emptyset$ iff
$X^{G}\neq \emptyset$.}
\vspace{.2in}

Let $G\subset\Gamma$ and let $L$ be the subcomplex of $K$ defined by\\  

$L=\{\sigma\in K\, :\,\, G\sigma =\sigma\}=\{\sigma \in K\, :\, \, 
G\,\,  fixes\,\,\sigma\,\, pointwise\}$ (by claim 1).\\

If $x\in K^{G}$  and $x\in int(\sigma )$, then $G\sigma =\sigma $ (the 
action is simplicial). Thus $K^{G}= L$.

Let ${\cal{W}}=\{U\cap X^{G}:\,\, U\in {\cal{U}},\,\, U\cap 
X^{G}\neq\emptyset \}$.  Recall that $X^{G}$ is convex.
Hence $\cal W$ is a cover of $X^{G}$ by convex compact subsets of 
$X^{G}$. Denote by $J$ the nerve of $\cal W$.
\vspace{.2in}

{\bf subclaim. } $L$ is homotopy equivalent to $J$.
\vspace{.2in}

Let $U$ be a vertex of $L$. Then $\gamma U=U$, for $\gamma\in G$. This 
means that
$ \gamma x=x$, for $\gamma\in G$, where $x$ is the center of the ball 
$U$. Thus
$ U\cap X^{G}\neq\emptyset$. Conversely, if $ U\cap X^{G}\neq\emptyset$, then 
$\gamma U\cap U\neq \emptyset$, hence $\gamma U=U$, for $\gamma\in G$. This 
defines a surjection $U\mapsto  U\cap X^{G}$, from the vertices of $L$ 
to the vertices of $J$.

Now, if $\{U_{0},...U_{k}\}$ forms a simplex in $L$  (that is, if $\bigcap 
U_{i}\neq\emptyset$)
then $G[\bigcap U_{i}]=\bigcap U_{i}$.
 But then $Gp=p$, where $p\in \bigcap U_{i}$ is the center of the 
compact convex set $\bigcap U_{i}$  (the center is unique, see \cite{BGS}, 
p.10). Thus  $\bigcap (U_{i}\cap X^{G})=\bigcap U_{i}\cap 
X^{G}\neq\emptyset$. 

Hence, the surjection $U\mapsto  U\cap X^{G}$, from the vertices of $L$ 
to the vertices of $J$ defines a simplicial map $h: L\rightarrow J$. We 
show that
$h^{-1}(\sigma)$ is a simplex, for a simplex $\sigma$ of $J$. So, let 
$\sigma =\, <U_{1}\cap X^{G}, ...,U_{k}\cap X^{G}>\,\in J$. Then the vertices of 
$h^{-1}(\sigma)$ are $U'_{1},...,U'_{l}\in \cal{U}$, where, for every
$j=1,...,l$, there is some $i=1,...,k$ with $U'_{j}\cap X^{G}=U_{i}\cap 
X^{G}$.
But then $\emptyset\neq (\bigcap U_{i})\cap X^{G}=\bigcap U'_{j}\cap X^{G} 
\subset \bigcap U'_{j}$, which means that $h^{-1}(\sigma)$ is a simplex.

Consequently, $h$ is a proper cellular map. Hence $h$ is  a homotopy equivalence.
This proves the subclaim.
\vspace{.1in}

Note that the proof of the subclaim implies $K^{G}\neq \emptyset$ iff
$X^{G}\neq \emptyset$.
To finish the proof of the claim note that ${\cal{W}}$ is a brick 
decomposition of
$X^{G}$ in the sense of \cite{CK}: nonempty intersections of elements 
in ${\cal{W}}$ are compact and convex. Hence, by  lemma 1.1, they are AR's. Then the main 
result in \cite{CK} implies that $J$ is homotopy
equivalent to $X^{G}$ which is contractible. This completes the proof of the claim.
\vspace{.2in}

Let $\cal{C}$ be the collection of subgroups of $\Gamma$
given by 
$${\cal{C}}=\{G\,\,:\,\, X^{G}\neq\emptyset\,\}=
\{G\,\, :\,\, G\,\, subgroup\,\, of\,\, \Gamma^{x},\,\, for\,\, 
some\,\, x\in X\,\}$$

 The second statement of claim 2 implies
  
$${\cal{C}}=\{G\,\,:\,\, K^{G}\neq\emptyset\,\}=
\{G\,\, :\,\, G\,\, subgroup\,\, of\,\, \Gamma^{a},\,\, for\,\, 
some\,\, 
a\in K\,\}$$

 This, together with the first statement of claim 2, imply
that $K$ is a universal $(\Gamma , \cal{C}) $-complex.
Also note that, since $X^G$ is contractible and non-empty, for every 
$G\in \cal{C}$, $X$ is $\cal{C}$-free and $\cal{C}$-contractible. 
Hence  there is a $\Gamma$-map
$f: K\rightarrow X$ (see \cite{FJ}, p.286, and note that there it is 
not required for $X$ to
be a complex).

Recall also the definition of the canonical $g$ map of a space into its 
nerve: for $x\in X$,  let $U_{0},...,U_{k}\in{\cal{U}}$, be the set of 
all elements in $\cal U$ that contain
$x$. Then the barycentric coordinates $\lambda_{0},...,\lambda_{k}$
of $g(x)$ in the simplex $<U_{0},...,U_{k}>$ are
$$\lambda_{i}=\frac{d(x,X\setminus 
U_{i})}{\sum_{j=0}^{j=k}d(x,X\setminus U_{j})}$$
where $d$ denotes the metric on $X$.
Note that $g$ is a $\Gamma$-map: if $U_{0},...,U_{k}$, is the set of 
all elements in $\cal U$ that contain $x$, then $\gamma U_{0},...,\gamma 
U_{k}$, is the set of all elements in $\cal U$ that contain $\gamma x$, 
and
$$\frac{d(\gamma x,X\setminus \gamma 
U_{i})}{\sum_{j=0}^{j=k}d(x,X\setminus \gamma U_{j})}=\frac{d(x,X\setminus 
U_{i})}{\sum_{j=0}^{j=k}d(x,X\setminus U_{j})}$$

This means that $gf: K\rightarrow K$ is a $\Gamma$-map,  and from  $K$ 
being a universal $(\Gamma ,{\cal{C}})$-complex  we deduce  that $gf$ is
$\Gamma$-homotopic to the identity on $K$ (see \cite{FJ}, p.286, theorem A.2.)

We want to prove now that $fg$ is $\Gamma$-homotopic to the identity in 
$X$, but
we can not apply the argument above because we do not have that $X$ is 
a $\Gamma$-complex. But now we use the fact that $X$ is a 
CAT(0) space. 
\vspace{.2in}

{\bf Lemma 2.2.} {\it Let $Z$  be a $\Gamma$-space and  $X$ be a CAT(0) 
space on which $\Gamma $ acts by isometries. Then any two $\Gamma $-maps from 
$Z$ to $X$ are $\Gamma$-homotopic.}
\vspace{.2in}

{\bf Proof.} Let $f_{0},f_{1}: Z\rightarrow X$, be two $\Gamma$-maps.
Define $f_{t}(z)=\alpha_{z}(d(f_{0}(z),f_{1}(z))\, t)$, for $z\in Z$
and $t\in [0,1]$, where
$\alpha_{z}$ is the unique  geodesic beginning at $f_{0}(z)$ and 
ending at $f_{1}(z)$. It follows from well known facts about
geodesics on CAT(0) spaces (see for instance \cite{B}, \cite{BH}) that $f$ is  a continuous 
$\Gamma$-map.  
\vspace{.2in}

Since $X$ is proper, it is straightforward to show that all maps and 
homotopies are proper.
This completes the proof of the proposition.
\vspace{.2in}

{\bf Remark.} The condition of $X$ being proper is necessary. For take 
$X$ to be the cone over the integers and extend the action of the 
integers ${\Bbb{Z}}$ on itself to
the CAT(0) space $X$ with fixed point the vertex.  Then ${\Bbb{Z}}$ 
acts cocompactly by isometries and with discrete orbits on $X$ but $X$ 
is not {\it properly } $\Gamma$-homotopic to a $\Gamma$-finite 
$\Gamma$-simplicial complex. Note that $X$ is $\Gamma$-homotopic to a point, 
but not properly $\Gamma$-homotopic to a point.
\vspace{.3in}

Now, before applying the proposition to prove theorem C we need some 
comments and
a lemma.  Let $X$ be a proper CAT(0) space on which $\Gamma$ acts 
cocompactly by isometries and with discrete orbits and let $K$ be a 
$\Gamma$-finite $\Gamma$-simplicial complex
$\Gamma$-homotopy equivalent to $X$. Let $K_n$, $n=1,2,3...$ be a 
sequence of
subcomplexes of K such that they satisfy:
\vspace{.1in}

\hspace{1in}{\it (i) $K_{n+1}\subset K_{n}$}

\hspace{1in}{\it (ii) $\overline{K\setminus K_{n}}$ is compact}

\hspace{1in}{\it (iii) $\bigcap K_{n}=\phi$}
\vspace{.1in}

Let $\iota : K_{n+1}\ra K_{n}$ be the inclusion and denote by
$\cal{K}$ the
tower $\{ 
K_{1}\stackrel{\iota}{\leftarrow}K_{2}\stackrel{\iota}{\leftarrow}
K_{3}\stackrel{\iota}{\leftarrow}...\}$ (i.e. a inverse system 
indexed by
$\N$, see \cite{EG}). We say that $\cal K$ is a $(X,K)$-tower.
\vspace{.2in}

{\bf Lemma 2.3.} {\it Every $(X,K)$-tower is a polyhedral resolution of
$\partial X$.}
\vspace{.2in}

{\bf Remark.} For the definition of a resolution see \cite{M}. This  
lemma implies that the shape of $\partial X$ is determined by $\cal K$ in
{\it pro-H}, where {\it H } is the homotopy category.
\vspace{.2in}

{\bf Proof of the lemma.}
Let $d$ be the metric on X and fix a point $x_{0}\in X$. Write 
$B_{n}=B(x_{0},n)=
\{ x\in X\, : \, d(x,x_{0})\leq n\}$, $S_{n}=
\{ x\in X\, : \, d(x,x_{0})= n\}$ and $X_{n}=
\{ x\in X\, : \, d(x,x_{0})\geq n\}$, for $n=1,2,3...$
Let $\cal S$ be the tower $\{ 
S_{1}\stackrel{r}{\leftarrow}S_{2}\stackrel{r}{\leftarrow}
S_{3}\stackrel{r}{\leftarrow}...\}$, where the maps
$r:S_{n+1}\ra S_{n}$ are given by geodesic retraction. Because
$\partial X = lim_{\leftarrow}S_{n}$ and all $S_n$ are compact theorem 
5 of \cite{M}
implies that the tower ${\cal{S}}$ is a resolution of $\partial X$. But 
geodesic retraction also
induces an equivalence, in {\it pro-H }, between $\cal{S} $ and 
${\cal{X}}=\{ X_{1}\stackrel{j}{\leftarrow}X_{2}\stackrel{j}{\leftarrow}
X_{3}\stackrel{j}{\leftarrow}...\}$, where the $j$'s are the 
inclusions.
Hence $\cal X$ is also a resolution of $\partial X$. It remains to prove that
$\cal X$ and $\cal K$ are equivalent in {\it pro-H }. Let $f:X\ra K$ 
and $g:K\ra X$ be
$\Gamma$-maps such that $fg$ and $gf$ are 
$\Gamma$-homotopic to the corresponding identities. By 
taking subsequences we can assume that the following conditions hold.
\vspace{.1in}

\hspace{1in} {\it (1) $f(X_{n})\subset K_{n}$}

\hspace{1in} {\it (2) $g(K_{n+1})\subset X_{n}$}

\hspace{1in} {\it (3)  $gf\mid_{X_{n+1}}:X_{n+1}\ra X_{n}$ and 
$fg\mid_{K_{n+1}}:K_{n+1}\ra K_{n}$ are $\Gamma$-homotopic to 

\hspace{1.33in}the inclusions
$X_{n+1}\ra X_{n}$ and $K_{n+1}\ra K_{n}$, respectively.}
\vspace{.1in}

Hence, by (1) and (2) above we can define maps $F: {\cal{X}}\ra \cal K$ 
and 
$G: {\cal{K}}\ra \cal X$ given by $F_{n}=f\mid_{X_{n}}: X_{n}\ra K_{n}$ 
and 
$G_{n}=g\mid_{K_{n}}: K_{n}\ra X_{n-1}$. Also by (3) above we get 
that $GF:{\cal{X}}\ra \cal X$ is equivalent in {\it pro-H }
to the inclusion map ${\cal{X}}\ra\cal X$
(i.e. given by the inclusions $j:X_{n+1}\ra X_{n}$). But ${\cal{X}}=\{ 
X_{1}\stackrel{j}{\leftarrow}X_{2}\stackrel{j}{\leftarrow}
X_{3}\stackrel{j}{\leftarrow}...\}$, consequently the inclusion map 
${\cal{X}}\ra\cal X$
is equivalent in {\it pro-H } to the identity on $\cal X$. In the same 
way we also get that
$FG$ is equivalent to the identity on $\cal K$. This proves the lemma.
\vspace{.2in}

{\bf Proof of theorem C.} By proposition A, $X$ admits a $(X,K)$-tower 
$\cal K$. But $X$ is $\Gamma$-homotopy equivalent to $Y$, hence $Y$ is also
$\Gamma$-homotopy equivalent to $K$. All this together with lemma
2.3 imply that $\cal K$ is
a polyhedral resolution of both, $\partial X$ and $\partial Y$. This 
proves theorem C.
\vspace{.3in}

For a group $\Gamma$ acting on a space $Z$, let
${\cal{C}}_{Z}=\{G<\Gamma \,\,:\,\, Z^{G}\neq\emptyset\,\}$.
Corollary A follows from the following lemma.
\vspace{.2in}

{\bf Lemma 2.4.} {\it Let $X$ and $Y$ be proper CAT(0) spaces on which 
$\Gamma$ acts cocompactly by isometries and with discrete orbits.
Then the actions have the same isotropy 
(i.e. ${\cal{C}}_{X}={\cal{C}}_{Y}$)
if and only if $X$ is $\Gamma$-homotopy 
equivalent to $Y$.}
\vspace{.2in}

{\bf Proof.} 
Suppose that  $X$ is $\Gamma$-homotopy equivalent to $Y$.
Then there is a $\Gamma$-map $f:X\ra Y$. Let $g\in \Gamma$ and
$x\in X$ such that $gx=x$. Then $gy=y$, for $y=f(x)$. This implies
${\cal{C}}_{X}\subset {\cal{C}}_{Y}$. In the same way we obtain
${\cal{C}}_{Y}\subset {\cal{C}}_{X}$. Hence both actions have the same
isotropy.

Conversely, suppose ${\cal{C}}_{X}={\cal{C}}_{Y}={\cal{C}}$.
Then, by proposition A, $X$ and $Y$ are 
properly $\Gamma$-homotopy equivalent
to  $\Gamma$-finite $\Gamma$-simplicial complexes $K$ and $J$, 
respectively. In fact, in the proof of proposition A 
we showed that ${\cal{C}}_{X}={\cal{C}}_{K}$ and ${\cal{C}}_{Y}={\cal{C}}_{J}$.
Hence ${\cal{C}}_{K}={\cal{C}}_{J}={\cal{C}}$. But in the proof of proposition
A we also showed that $K$ and $J$ are universal $(\Gamma ,{\cal{C}})$-complexes,
which are unique, up to $\Gamma$-homotopy equivalence
(see \cite{FJ}, p.286). This proves the lemma
and corollary A.
\vspace{.2in}

{\bf Proof of Corollary B.} First recall that if $G$ is a finite group
acting on a proper CAT(0) space $X$, then $G$ fixes some point $p$.
This point is the unique center of the compact $G$-invariant set $Gp$
(see \cite{BGS}, p.10).
This implies $G\in {\cal{C}}_{X}$, for every $G$ finite.
Hence, if the isotropy ${\cal{C}}_{X}$ consists only of finite subgroups,
then it is the family of {\it all} finite subgroups.
Consequently, if $\Gamma$ acts on $X$ and $Y$ with finite isotropy
(i.e. all isotropy groups are finite), then ${\cal{C}}_{X}={\cal{C}}_{Y}=\{ G<\Gamma\, :\, G$ is finite $\}$.
The corollary now follows from corollary A.
\vspace{.5in}

{\bf 3. Proof of theorem A.}
\vspace{.2in}

{\bf Theorem A.} {\it Let $X$ be a noncompact proper CAT(0) space on 
which $\Gamma$ acts cocompactly
by isometries. If $H^{i}_{c}(X)\neq 0$, for some $i$,
then $X$ is almost geodesically complete.}
\vspace{.3in}

{\bf Proof.}
Here is the idea of the proof. If $X$ is not almost geodesically complete we 
construct (see claim 2 below) retractions $f_{r},\,\, r>0$, properly homotopic
to the identity, such that $f_{r}(X)$ misses an arbitrarily large ball about a fixed point $p$
(the balls get larger as $r\ra\infty$).
Then we prove in claim 3 that $H^{*}_{c}(X) = 0$, in the following way. 
First we prove that for
$z\in H^{*}_{c}(X)$ there is an $r>0$ such that $f_{r}^{*}(c)=0$
(for this just take $r$ large enough so that $f_{r}(X)$ misses the 
support of $z$). But $f_{r}$ is properly homotopic to  the identity
$1_{X}$, hence $z=1_{X}^{*}z=f_{r}^{*}z=0$.
We now give a detailed proof.

Given $p\in X$ and $s>0$ define $f^{p,s}: X\rightarrow X$, by
$$f^{p,s}(\alpha (t))=\left\{\begin{array}{ll}p=\alpha (0) & t\leq s \\ 
\alpha (t-s) &
t\geq s\end{array}\right.$$
where $\alpha$ is a geodesic  beginning at $p$.
(This makes sense because for every two points in a CAT(0) space there
is a unique geodesic segment joining these two points and depending 
continuously on them. See \cite{B}, \cite{BH}.)
Note that $f^{p,s}$ is a proper map properly homotopic to the identity 
in $X$. 

Also, given a geodesic segment $\alpha$ denote by $\ell_{\alpha}$ the 
supremum
over all $\ell$ such that $\alpha$ can be extended to the interval 
$[0,\ell ]$.
If the maximum does not exist we write $\ell_{\alpha}=\infty$. Note 
that
if $\ell_{\alpha}<\infty$, there is a geodesic segment defined on the 
interval $[0,\ell_{\alpha}]$ 
extending $\alpha$  (because $X$ is proper), and if $\ell_{\alpha}=\infty$ then $\alpha$  
can be extended to $[0,\infty)$. (For this last statement we can use a 
general Arzel\'a Ascoli argument,
or, for simplicity in this special case, the fact that  $X\cup\partial X $
is compact. For the definition of the boundary $\partial X$ and its 
properties, see for instance, \cite{BH}).

 Fix $p\in X$. Now, since the action is cocompact, to prove the 
theorem it is enough to prove the following:
\vspace{.2in}

{\it  There is a  constant $C$ such that for every $\gamma\in \Gamma$ 
there is a geodesic ray $\alpha :
[0,\infty)\rightarrow X$, $\alpha (0)=p$, and $d(\gamma (p),\alpha)\leq 
C$.}
\vspace{.2in}

We prove this by contradiction. Assume that the statement above
does not hold for $X$. We have 
\vspace{.2in}

(*) Given $r>0$ there is $\gamma_{r}\in \Gamma$ such that 
$\ell_{[p,x]}<\infty$ for $x\in B_{\gamma_{r}(p)}(r)$. 
\vspace{.2in}

Here $B_{\gamma_{r}(p)}(r)$ denotes the closed ball with
center $\gamma_{r}(p)$ and radius $r$ and $[p,x]$ denotes the (unique) geodesic 
segment from $p$ to $x$.
\vspace{.2in}

{\bf Claim 1.} Given $r>0$,  we have $$sup\, \{\ell_{[p,x]}\,\,   : 
\,\, x\in B_{\gamma_{r}(p)}(r)\}\,\,<\,\, \infty$$

We prove this by contradiction. Suppose there is a sequence 
$\{x_{j}\}\subset B_{\gamma_{r}(p)}(r)$ with $\ell_{[p,x_{j}]}\rightarrow\infty$ 
and consider the sequence
$\{y_{j}\}$, defined by $y_{j}=\alpha (\ell_{[p,x_{j}]})$, where 
$\alpha$ is a geodesic segment extending $[p,x_{j}]$ (choose one $\alpha$). 
Because $X\cup\partial X$ is compact we can assume that the sequence $\{y_{j}\}$ converges to a point 
$x_{0}$ in the boundary $\partial X$ of $X$. Then the geodesic 
ray $[p,x_{0})$ intersects $B_{\gamma_{r}(p)}(r)$,  which contradicts 
(*).
\vspace{.3in}

{\bf Claim 2.} Given $r>0$ there is $f_{r}:X\rightarrow X$,  properly 
homotopic to the identity, such that $f_{r}(X)\subset X\setminus B_{p}(r)$ 
\vspace{.2in}

Note that, since $X$ is noncompact (and claim 1), we can assume that 
$p\notin B_{\gamma_{r}(p)}(r)$.
Write   $s=sup\, \{\ell_{[p,x]}\,\,   : \,\, x\in 
B_{\gamma_{r}(p)}(r)\}\,\,<\,\, \infty$,
and remark that $f^{p,s}(X)\subset X\setminus B_{\gamma_{r}(p)}(r)$. 
Hence $\gamma_{r}^{-1}f^{p,s}(X)\subset X\setminus B_{p}(r)$. Take 
$f_{r}=\gamma_{r}^{-1} f^{p,s}$.
This proves the claim.
\vspace{.2in}

{\bf Claim 3.} $H^{i}_{c}(X)=0$, for all $i$.
\vspace{.2in}

Take a $z\in H^{i}_{c}(X)$. Then there is a $r>0$ such that
$\iota (z)=0$, where $\iota : H^{i}_{c}(X)\rightarrow 
H^{i}_{c}(X\setminus int\, B_{p}(r))$
 is the restriction map. Take $f_{r}$ from claim 2 and consider the 
following
commutative diagram corresponding to the map of pairs 
$f_{r}:(X,X)\rightarrow (X,X\setminus int\, B_{p}(r))$,

\begin{picture}(250,150)(-100,0)
\put (0,40){$H_{c}^{i}(X)$} 
\put (0,120){$H_{c}^{i}(X)$}
 \put (120,120){$H_{c}^{i}(X\setminus int\, B_{p}(r))$}   
 \put (120,40){$H_{c}^{i}(X)$}
\put (0,90){$f_{r}^{*}$}
 \put (120,90){$f_{r}^{*}$} 
\put (20,110){\vector(0,-1){50}}
\put (140,110){\vector(0,-1){50}}
   \put (40,125){\vector(1,0){70}}
\put (40,45){\vector(1,0){70}}
\put (70,130){$\iota$}
\put (70,50){$1_{X}$}
\end{picture}

\noindent where $1_{X}$ denotes the identity. Then 
$f_{r}^{*}(z)=f_{r}^{*}(\iota(z))=0$. But
$f_{r}$ is properly homotopic to the identity, thus $f_{r}^{*}$ is the
identity
and $z=0$. This completes the proof of Theorem A.
\vspace{1in}

{\bf 4. Proof of Proposition B.}
\vspace{.2in}

Before proving theorem B, we recall some definitions and results
from PL topology (see \cite{RS}).

Let $X$ be a $PL$ space. A triangulation on $X$ is a pair $(T,\rho )$, 
where $T$ is a simplicial
complex and $\rho : T\rightarrow X$ is a $PL$ homeomorphism. Sometimes we 
just say that $T$ is a 
triangulation of $X$. For a simplicial complex $T$, $\mid T\mid$ 
denotes the underlying topological space. From now on, all simplicial 
complexes are locally finite and finite
dimensional.

Let $K\subset T$ be simplicial complexes. We say that $K$ is full in 
$T$ if every simplex of $T$ intersects $K$ in exactly one (possibly empty) face.
Define $N(T,K)$ to be the subcomplex of $T$ that contains all simplices
intersecting $\mid K \mid$, together with their faces and $C(T,K)$ the 
subcomplex of
 $T$ that contains all simplices not
intersecting $\mid K \mid$. Define also $int \, N(T,K)=\, \mid 
N(T,K)\mid \, \setminus
\,  \mid C(T,K)\mid $.

Let $Y$ be a $PL$ space and $X\subset Y$ a $PL$ subspace.
 Then $N=\,\, \mid N(T',K) \mid$ is called a regular neighborhood
of $\mid K\mid$ in $\mid T\mid$, where: $T$ is a triangulation of $Y$ 
and
$K$ the subcomplex of $T$ triangulating $X$ with $K$  full in $T$,
and $T'$ is a first derived of $T$ near $K$
(i.e $T'$ is obtained from $T$ by subdividing only simplices that are 
neither in $K$,
nor in $C(T,K)$, see p.32 of \cite{RS}).

Let  $\Delta^{n}$ and $\Delta^{m}$ be  standard simplices. Then 
$\Delta^{n}*\Delta^{m}=\Delta^{n+m+1}$, where the star denotes ``join" (see 
\cite{RS}, p.2).
Note that  every 
$z\in\Delta^{n+m+1}$ can
be written uniquely as $z=(1-s)x+sy$, $s\in [0,1]$, $x\in\Delta^{n}$, 
$y\in\Delta^{m}$ (here we are assuming $\Delta^{n}, 
\Delta^{m}\subset\R^{N}$ in general position). Note also that for $s=0$ 
we get  the points of $\Delta^{n}$, and  for $s=1$ 
we get  the points of $\Delta^{m}$.

Define, for $t\in [0,1]$, the canonical deformation retraction of $c_{t}: 
(\Delta^{n+m+1}\setminus \Delta^{m})\rightarrow
(\Delta^{n+m+1}\setminus \Delta^{m})$ by $c_{t}((1-s)x+sy)=((1-ts)x+tsy)$,
 $x\in\Delta^{n}$, 
$y\in\Delta^{m}$, $s\in [0,1)$ (note that $c_{t}$ is well defined and 
continuous for $s\neq 1$). 
Hence $c_{1}$ is the
identity, $c_{0}(\Delta^{n+m+1}\setminus \Delta^{m})\subset \Delta^{n}$ 
and $c_{t}(x)=x$ for all
$t\in [0,1]$ and $x\in\Delta^{n}$. 
\vspace{.1in}

{\bf Remarks.}

{\bf 1.} Let $d$ be any metric on the simplex
$\Delta^{n+m+1}$, compatible with the topology of 
$\Delta^{n+m+1}$. We will use the following simple fact:
given $\delta >0$ there is an open neighborhood 
$V_{\delta}$ of $\Delta^{n}$ in 
$\Delta^{n+m+1}$, such that $d(c_{1}(x),c_{1}(y))\leq d(x,y)+\delta$, 
for $x,y \in V_{\delta}$.
(Proof: because $c_{1}(x)=x$ for $x\in \Delta^{n}$ there is a 
neighborhood 
$V_{\delta}$ of $\Delta^{n}$ such that $d(c_{1}(x),x)\leq \delta /2$, 
for $x\in V_{\delta}$.
Then, for $x,y\in V_{\delta}$ we get $d(c_{1}(x),c_{1}(y))\leq 
d(c_{1}(x),x)+ d(x,y)+
d(y,c_{1}(y))\leq \delta + d(x,y)$.)
\vspace{.05in}

{\bf 2.} Let $K$ be full in $T$. Then every simplex $\Delta$ in $N(T,K)\setminus 
(C(T,K)\cup K) $ can be 
written uniquely
as $\sigma*\tau$, where  $\sigma< \Delta$ is the simplex
$\Delta \cap K$ and $\tau$ is the {\it complementary} simplex $\Delta\cap 
C(T,K)$. Hence  using $c_{t}$ defined
above, we can construct simplexwise a deformation retraction 
$c_{t}=c_{t}(T,K): int (N(T,K))
\rightarrow  int (N(T,K))$ with the same properties, that is, $c_{1}$ 
is the
identity, $c_{0}(int (N(T,K)))\subset\, \mid K\mid$ and $c_{t}(x)=x$ 
for all
$t\in [0,1]$ and $x\in | K|$. Write $c=c_{0}$.
\vspace{.2in}

Let $T$ be a triangulation of the $PL$ space $Y$. For every simplex 
$\sigma\in T$ choose a flat  metric $d_{\sigma}$ on it (i.e. there is a 
linear isometric embedding from $(\sigma,d_{\sigma})$ into  some 
euclidean space) in such a way that the metrics  coincide on each intersection 
of simplices. This gives a way to define the length of a path in $Y$ 
and this determines a metric $d$  on $Y$ by defining
$d(x,y)=inf\{$lengths of paths joining $x$ to $y$ $\}$. We say that $d$ 
is a piecewise flat metric on the $PL$ space $Y$ (or a $PL$ metric on 
$Y$), with respect to the triangulation $T$.
We have that if $T$ is locally finite and ($Y$,$d$) is complete
(as a metric space) then $Y$ with metric $d$ is a geodesic space.
(Recall that we are assuming all simplicial complexes to be locally 
finite.)
\vspace{.12in}

Let $Y$ be a $PL$ space with a  $PL$ metric $d$, where $d$ is piecewise 
flat with respect
to the triangulation $T$ of $Y$, and let $L$ be a subcomplex of $T$. We 
denote by $d_{L}$
the intrinsic metric on $\mid L\mid$ induced by $d$, i.e. 
$d_{L}(x,y)=inf\{$lengths of paths in L joining $x$ to $y$ $\}$.
Let $U$ be any other triangulation of $Y$ and $J$ a subcomplex of $U$, 
with $|J|\subset
| L|$. Define $mesh_{L}(J)=sup\{ diam_{L}(\sigma )\,\, 
:\,\,\sigma\in J\}$, where
$diam_{L}(\sigma )$ is the diameter of $|\sigma|\subset| L|$ with respect to the metric $d_{L}$. 
If $L=T$ we simply write $mesh$ instead of $mesh_{L}$.
Also, define $mesh_{0}(J)=sup\{ length(\sigma )\,\, :\,\, \sigma\in 
J^{1}\}$, where
$J^{1}$ is the set of one simplices in $J$. Note that if $J,L$, $|J| \subset
| L|$, are subcomplexes of
a subdivision of $T$, then $mesh_{L}(J)\leq mesh_{0}(J)\leq 
mesh_{0}(T)$.
\vspace{.2in}

{\bf 4.1. Lemma.} {\it Let $N_{1}$ and $N_{2}$ be regular neighborhoods 
of the $PL$ space
$X$ in the $PL$ space $Y$, with piecewise flat metric $d$. Then there is 
a $PL$ homeomorphism
$j: Y\rightarrow Y$, with $j(N_{1})=N_{2}$, which is the identity
on $X$. Moreover
if \newline $N_{i}=\,\,| N(T_{i}',K_{i})|$, where $K_{i}$ is 
induced by $T_{i}$, $i=1,2$,  we have $$d(x,j(x))\leq 
mesh_{Y}(T_{1}')+mesh_{Y}(T_{2}')$$}
\vspace{.1in}
(Here $T_{i}'$ is a first derived of $T_{i}$ near $K_{i}$, $i=1,2$.)
\vspace{.2in}

{\bf Proof.} The first part is theorem 3.8 of \cite{RS} (drop 
compactness). The inequality follows because
the map constructed in the proof of theorem 3.8 of \cite{RS} is a 
composition of two maps,
each moving a point in $\sigma\cap N_{i}$, where  $\sigma $ is a simplex of 
$N(T_{i},K_{i})$, $i=1,2.$ This proves lemma 4.1.
\vspace{.2in}

Denote the map $j$ given by the lemma as $j(T_{1}',T_{2}')$.
\vspace{.2in}

We will also need the following result (see \cite{RS}, p.32)
\vspace{.2in}

{\bf 4.2. Lemma.} {\it Let $X$ be a $PL$ subspace of the $PL$ space 
$Y$,
$T$ a triangulation of $Y$ and $K$ a subcomplex of $T$ triangulating 
$X$,
where $K$ is full in $T$. Given an open neighborhood $U$ of $X$ in $Y$
there is a first derived $T'$ of $T$ near $K$, such that the regular 
neighborhood
$| N(T',K)|$ is contained in $U$.}
\vspace{.2in}

{\bf 4.3. Proposition B.} {\it  Let K be a locally finite, 
contractible
simplicial complex which admits a cocompact simplicial action. Then 
$H^{i}_{c}(\, | K|\,  )\neq 0$, for some $i$.}
\vspace{.2in}

{\bf Proof.} 
First we give some motivation and the idea of the proof.

We will prove the proposition by contradiction, so we will assume 
$H_{c}^{i}(\, | K| \,
)=0$, for all $i$. Now, consider first a particular case:
suppose that $K$ is $PL$ homeomorphic to a $PL$ manifold $M^{n}$. 
Since $M^{n}$ is contractible and $H_{c}^{i}(M)=0$, for all $i$,
we can assume, after some stabilization,
that $M^{n}$ is $PL$-homeomorphic to euclidean half $n$-space
$\R^{n}_{+}=\{ (x_{1},...,x_{n})\,\, :\,\, 
x_{n}\geq 0 \}$ and $\partial M$ is $PL$-homeomorphic to $\R^{n-1}$.
Denote by $\Gamma$ the group that acts
simplicially and cocompactly on $K$ and let $d$ be a $\Gamma$-invariant 
metric on $K\cong_{PL}M$. We will prove these two key facts:

(a) $M$ is $d$-close to $\partial M$ (this is {\bf (2)} in the proof),

(b) bounded extensions of maps to $\partial M$(for a precise statement see claim 4).
\vspace{.05in}

Then we proceed as follows (see arguments after the proof of claim 4).
Take $x_{0}\in \partial M$. Let $D$ denote a closed $PL$ $n$-ball. We can construct a $PL$ embedding
$\rho : (D , \partial D )\ra (M,\partial M )$ such that

(i) $\rho (D )$ is far from $x_{0}$

(ii) $\rho|_{\partial D}\neq 0\in \pi_{n-2}(\partial M\setminus\{ 
x_{0}\} )$.
\vspace{.05in}

Using (a), (b) and $\rho$ above we construct a map $\psi :D\ra \partial M$ with 
$\rho|_{\partial D}=\psi|_{\partial D}$ and $\psi 
(D)\subset \partial M\setminus\{ x_{0}\}$. Hence 
$\rho|_{\partial D}=\psi|_{\partial D}=0\in
\pi_{n-2}(\partial M\setminus\{ 
x_{0}\} )$, a contradiction.

For the general case, $K$ may not be $PL$ homeomorphic to a $PL$ manifold.
Hence we ``replace" $K$ by its regular neighborhood $M$, in some euclidean 
space. We do this at the beginning of the proof.
The problem now is that $\Gamma$ does not act, at least a priori, 
simplicially on $M$. To prove
(b) above (i.e. claim 4 in the proof) we use an ``approximate" action of $\Gamma$ on $M$
(see claim 2 in the proof).
\vspace{.1in}

We now give a detailed proof.
\vspace{.1in}

The proof is by contradiction. So suppose  
$H^{i}_{c}(\, | K|\, )= 0$, for all $i$, and denote by $\Gamma$ the group that acts
simplicially and cocompactly on $K$.

Now, note that, since $K$ is locally finite and connected, $K$ is 
countable, that is, it has a countable number of simplices. Note also that
$K$ is finite dimensional. Hence we can
embed $K$ simplicially and properly
in some ${\Bbb{R}}^{n}$, $n\geq 2(dim\, K)+2$, and take $n\geq 6$.
Let $T$ be a triangulation of ${\Bbb{R}}^{n}$ such that there is 
a full subcomplex $J$ of $T$
with $| K|\,\, =\,\, | J|$ and  $J$ is a 
subdivision of $K$. Give $K$ the unit metric, that is, the geodesic 
piecewise flat metric
where every edge has length one. Note  that this metric is   
$\Gamma$-invariant and that $K$, with this metric, is proper.
Denote this piecewise flat metric on $| J|\, =\, | K|$
by $d_{K}$. By lemma 2.5 of \cite{BHu}
(the proof works for infinite  complexes), maybe after a subdivision 
of $T$ away from $K$, we can extend this proper piecewise flat metric 
to a proper piecewise flat metric $d$ on ${\Bbb{R}}^{n}=| T|$.  Also, maybe 
after further subdivision, we can assume $mesh_{0}(T)\leq 1$, and that
every simplex of $J$ is convex in $| T|$ (see \cite{BH}).

Recall that $T'$ is a first derived of $T$ near $J$ (see \cite{RS}, p.32).
Let $M=\, | N(T',J)|$ be a regular neighborhood of $| J| $ 
in  ${\Bbb{R}}^{n}$.
Then $M$ is a $n$-manifold with boundary $\partial M$ and $M$ is 
properly homotopic to $| K| $. Hence $M$ is contractible and 
$H^{i}_{c}(M)=H^{i}_{c}(\, | K| \, )=0$, for all $i$. By duality we have that 
$H_{i}(\partial M)=0,\,\, i>0$.  Note
that because $M$ is contractible, the boundary of the regular 
neighborhood of $M$ in ${\Bbb{R}}^{n+1}$ is the suspension of the regular 
neighborhood of $M$ in ${\Bbb{R}}^{n}$,
where the embedding $M\ra\R^{n+1}$ is the composition 
$M\ra\R^{n}\ra\R^{n+1}$. Hence, by embedding canonically ${\Bbb{R}}^{n}$ into 
${\Bbb{R}}^{n+1}$ if necessary, we can assume that $\partial M$ is 
simply connected. This implies, together with $H_{i}(\partial M)=0,\,\, i>0$,
that $\partial M$ is contractible. 

We claim that we can also assume that $M$ and $\partial M$ are simply connected at 
infinity. To see this, just cross $\Gamma$  with 
${\Bbb{Z}}^{2}$ and $K$ with ${\Bbb{R}}^{2}$ if necessary 
and make $\Gamma\times \Z^{2}$ act cocompactly on $K\times\R^{2}$. Note 
that we still have $H_{c}^{i}(K\times\R^{2})=0$, for all $i$.
Note also that $X\times \R^{2}$ is simply connected at infinity for $X$ 
simply connected.

Recall that if $X^{m}$ is a contractible simply-connected-at-infinity 
high-dimensional $PL$-manifold with empty boundary, then $X$ is $PL$-homeomorphic to 
euclidean $m$-space.
Moreover, if $X^{m+1}$ is a contractible simply-connected-at-infinity 
high-dimensional $PL$-manifold with boundary $PL$-homeomorphic to euclidean $m$-space,
then $X$ is $PL$-homeomorphic to half euclidean $(m+1)$-space (see 
\cite{S}, \cite{BLL}). It follows from these last remarks that
we can assume $\partial M$ to be $PL$-homeomorphic to 
euclidean $(n-1)$-space
${\Bbb{R}}^{n-1}$ and $M$ to be $PL$-homeomorphic to euclidean half $n$-space
${\Bbb{R}}^{n}_{+}=\{x=(x_{1},...,x_{n})\in {\Bbb{R}}^{n}\,:\,\, x_{1}\geq 0\}$.

Denote by $d_{M}$ the  intrinsic metric on $M$ induced from 
${\Bbb{R}}^{n}$ with 
 metric $d$.
Because $|J|\subset M$ we have $d_{M}|_{J}\leq d_{K}$. 

Let $c$ be the retract $c=c(T,J): int (N(T,J))
\rightarrow  int (N(T,J))$ defined in the introduction of this section.
Note that $c=c(T,J)$ depends on the choice of $M= | N(T',J)|$,
and the choice of $M$ depends only on the choice of the first derived $T'$.

Remark that $(M,d_{M})$ is a proper metric space. (This is because
$T$ is finite dimensional locally finite and proper, hence the same is
true for any subcomplex of a  (locally finite)  subdivision of $T$.)

We prove now that we can choose $M$ close to $| J|$ to get $d_{K}$
close to $d_{M}|_{J}$. In fact $M$ will get closer and closer to
$| J|$, as we approach infinity.
\vspace{.2in}

{\bf Claim 1.} We can choose $M$   so 
that $d_{K}\leq (d_{M}|_{J})+1$
\vspace{.1in}

Enumerate all simplices of $J$: 
$\sigma_{1},\sigma_{2},...$, and let $\Delta_{1}, \Delta_{2},...$ their 
corresponding simplices in $N(T,J)$ of maximal dimension (see remark 2 at the beginning 
of this section).
Let $A=A(n)$ be the number of simplices of the first barycentric 
subdivision of a $n$-simplex $\Delta^{n}$.
Remark 1 at the beginning 
of this section imply that
there  are open neighborhoods 
$V_{i}$ of $\sigma_{j}$ in $\Delta_{j}$, such that for all 
$x,y\in V_{i}$ we have
$$d_{\sigma_{j}}(c(x),c(y))\leq  d_{\Delta_{j}}(x,y)+\frac{1}{A2^{j}}$$

Define $W=int\,\bigcup V_{i}$. Since all triangulations here are locally 
finite, we have that $W$ is an open neighborhood
of $| J|$. We can assume $\overline{W}\subset int\, N(T,J)$, so that $c$ is
defined at all points of $W$.

By lemma 4.2 we can choose a first derived $T'$ of $T$ near $J$ such 
that $M=\, | N(T',J)|\subset W$. Note that $c$ is defined
for every point in $M\subset W$.
From the definition of first derived follows that for every $\Delta_{j}\in 
N(T,J)$ there is at most $A$ simplices $\Delta '\in N(T',J)$ with 
$\Delta '\subset \Delta_{j}$. 

Since we took $T$ small enough so that every simplex is convex,
we have that $d_{K}(x,y)=d_{\sigma_{j}}(x,y)$ if $x,y\in\sigma_{j}\in J$
and $d_{M}(x,y)=d_{\Delta '}(x,y)$ if $x,y\in\Delta '\in N(T',J)$. Hence we can
write the above inequality in the following form:

$$d_{K}(c(x),c(y))\leq  d_{M}(x,y)+\frac{1}{A2^{j}}$$

\noindent for $x,y\in \Delta '\in N(T',J)$ and $\Delta '\subset\Delta_{j}\in N(T,J)$.

Now, for $x,y\in \,| J| $, let $\alpha : [0,a]\rightarrow M$,
$\alpha (0)=x$ and $\alpha (a)=y$, be a distance minimizing path
with respect to $d_{M}$, and  let $0=t_{0}<t_{1}< ...<t_{r}=a$ be such 
that for each
$i=1,...,r$ there is a simplex $\Delta '_{i}\in N(T',K)$,
 with $\alpha (t_{i-1}),\alpha (t_{i})\in\Delta_{i}'$ and all 
$\Delta_{i}'$
different. Note that we can choose all $\Delta_{i}'$
different because  every simplex of $T'$ is convex.
Also, for each $i$, let $j_{i}$ be such that $\Delta_{i} '\subset\Delta_{j_{i}}$
Write $x_{i}=\alpha (t_{i})$. 
Thus $d_{M}(x,y)=\sum d_{M}(x_{i-1},x_{i}).$ For $x,y\in \,| J| $ we have
\vspace{.1in}

$$ d_{K}(x,y)\leq \sum d_{K}(c(x_{i-1}),c(x_{i}))
\leq 
 \sum d_{M}(x_{i-1},x_{i})+\frac{1}{A2^{j_{i}}} \leq d_{M}(x,y)+1$$

For the first inequality use triangular inequality plus the fact that
$c(x_{0})=c(x)=x$ and $c(x_{r})=c(y)=y$. For the last inequality use the 
fact that for each $j_{i}$ there are, at most, $A$ simplices $\Delta_{i} '\in N(T',J)$ with 
$\Delta_{i} '\subset \Delta_{j_{i}}$.

This proves claim 1.
\vspace{.2in}

Note that, because we can make $ mesh_{0}(T)$ small,  we can assume
\vspace{.2in}

{\bf (1)}\hspace{.4in} $d_{M}(x,c(x))\leq 1$,  and
\vspace{.2in}

{\bf (2)}\hspace{.2in}  for every $x\in M$ there is $w\in\partial M$ 
with
$d_{M}(x,w) \leq 1$ (take $w=c(x)$ and use {\bf (1)}).
\vspace{.1in}

Recall that $\Gamma$ acts on $(K,d_{K})$ and on $(J,d_{K})$ simplicially by isometries.
\vspace{.2in}

{\bf Claim 2.} For every $\gamma\in \Gamma$, there is a $PL$ 
homeomorphism
$g: M\rightarrow M$, such that $d_{K}(cg(x),\gamma 
c(x))\leq 15$, for all $x\in M$.
\vspace{.2in}

Let $\iota :| J| \rightarrow {\Bbb{R}}^{n}$ denote the inclusion.
Because we assumed $n\geq 2(dim\, J)+2$ the two embeddings $\iota$ and 
$\iota\gamma$ are ambient isotopic. Then, by the uniqueness of regular 
neighborhoods,   
 there is a $PL$-homeomorphism $g_{1}:M\rightarrow M$, such that 
$g_{1}|_{K}=\gamma$.

Now, because
$cg_{1}c_{t}(x)\rightarrow cg_{1}c(x)=\gamma c(x)$, as $t\rightarrow 
1$,
for every $x\in M$, there is an open neighborhood $U\subset M$ of $| 
J|$ such that,
\vspace{.2in}

${\bf (3)}\hspace{1in} d_{K}(cg_{1}(x),\gamma c (x))<1\,\,\,\,\,\,\, 
x\in U $
\vspace{.2in}

By lemma 4.2 there is a first derived $T''$ of $T$ near $J$ with
\vspace{.1in}

${\bf (4)}\hspace{1in} N\, =\,| N(T'',J)|\,\,\subset\, U$.
\vspace{.1in}

Note that, because $mesh_{0}(T)\leq 1$, we have that $mesh_{M}(T')\leq 
1$ and
$mesh_{M}(T'')\leq 1$ (see comment right before lemma 4.1). 
Now, let $h=j(T',T'')$, given by lemma 4.1, and we get 
$d_{M}(x,h(x))\leq 2$,
for $x\in M$.
This, together with {\bf (1)} and claim 1, imply, for $x\in M$,
\vspace{.1in}

${\bf (5)}\hspace{.66in}d_{K}(ch(x),c(x))\leq  d_{M}(ch(x),c(x))+1$

$\hspace{2in}\leq d_{M}(ch(x),h(x))+d_{M}(h(x),x)+d_{M}(x,c(x))+1\leq 5$
\vspace{.1in}

Let also $j=j(g_{1}T'', T')$, so $j(g_{1}N)=M$. If $\sigma $ is a 
simplex in $N(T'',J)$,
 {\bf (1)},{\bf (3)},
{\bf (4)} imply:
\vspace{.1in}

$\hspace{1.08in}diam_{M} ( g_{1} \sigma ) 
\leq 
diam_{M}( cg_{1}\sigma   ) + 2$
\vspace{.1in}

$\hspace{1.99in}\leq diam_{K}(   cg_{1}\sigma   ) + 2
\leq diam_{K}( \gamma c   \sigma   ) +4$
\vspace{.1in}

$\hspace{1.99in}\leq diam_{K}(  c   \sigma   ) + 4 \leq mesh_{0}(T) + 4 
\leq 5.$
\vspace{.1in}

Here $diam_{K}(A)$ is the diameter of $A\subset| J| =| K|$, 
with metric $d_{K}$.
\vspace{.1in}

Hence, by lemma 4.1 we get 
\vspace{.1in}

$d_{M}(x,j(x))\leq mesh_{M}(g_{1}T'')+mesh_{M}(T')\leq 5+ 1= 6$
\vspace{.1in}

This together with {\bf (1)} and claim 1 imply for $x\in g_{1}N$,
\vspace{.1in}

${\bf (6)}\hspace{.5in} d_{K}(c(x), cj(x))\leq \,\,\,${\Large [}$\,\, d_{M}( 
c(x) ,x  ) + d_{M}(x  ,j(x)  ) 
+ d_{M}(j(x)  ,cj(x)  )\,\,  ${\Large ]} \, +\, 1
\vspace{.1in}

$\hspace{1.82in}\leq (1+6+1)+1=  9$
\vspace{.1in}

Define $g=jg_{1}h$. Then {\bf (3)},{\bf (4)},{\bf (5)}, and  {\bf (6)} 
imply, for $x\in M$,
\vspace{.1in}

$d_{K}( cg(x),\gamma c(x)  )= d_{K}( cjg_{1}h(x),\gamma c(x)  )$
\vspace{.1in}

$\hspace{1.2in}\leq
 d_{K}(cj[g_{1}h(x)] , c[g_{1}h(x)] )+
 d_{K}(cg_{1}h(x) , \gamma c(x)  )$
\vspace{.1in}

$\hspace{1.2
in}\leq  9+ d_{K}(cg_{1}[h(x)] , \gamma c[h(x)]  )+d_{K}(\gamma 
ch(x),\gamma c(x))$
\vspace{.1in}

$\hspace{1.2in}\leq 9+1+
d_{K}(\gamma ch(x),\gamma c(x))=10+d_{K}(ch(x),c(x))\leq 10+5=15$.
\vspace{.1in}

 This proves claim 2.
\vspace{.2in}

{\bf Claim 3.} Given $b>0$,  there is an $a_{b}$, such that the 
following holds.
For any map $\phi :\partial \sigma\rightarrow \partial M$, where 
$\sigma$ is a simplex,
 with $diam_{M}(\phi (\partial \sigma))\leq b$, there is an extension 
$\phi : \sigma\rightarrow \partial M$  with $diam_{M}(\phi ( 
\sigma))\leq a_{b}$.
\vspace{.2in}

Let $C_{0}$ be a finite subcomplex of $J$ such that for any subset $C'$ 
of $J$, with
$diam_{K}(C')\leq b+3$, there is a $\gamma\in\Gamma$ such that 
$C'\subset \gamma C_{0}$. Let also $C_{1}$ be a finite subcomplex of $J$ such 
that $d_{K}(J\setminus C_{1},C_{0})\geq 16$.

Consider now $c^{-1}(C_{1})\cap\partial M$. There is a $B\subset 
\partial M$, homeomorphic
to a $(n-1)$ ball such that  $\partial M \cap c^{-1}(C_{1})\subset B$ 
(recall that $\partial M$ is homeomorphic to
${\Bbb{R}}^{n-1}$ and $c$ is proper). Let $E_{0}$ be a finite subcomplex of $K$ such that 
$B\subset c^{-1}(E_{0})$. Then for any map $\partial\sigma\rightarrow 
c^{-1}(C_{1})\cap\partial M\subset B$ there is an extension $\sigma\rightarrow 
c^{-1}(E_{0})\cap\partial M$. Let also $E_{1}$ be a finite subcomplex of 
$K$ such that $d_{K}(K\setminus E_{1},E_{0})\geq 16$.

Let $\phi:\partial\sigma\rightarrow\partial M$, with 
$diam_{M}(\phi(\partial\sigma ))\leq b$. Then (use claim 1 and {\bf 
(1)})$$diam_{K}(c\phi(\partial\sigma ))\leq  diam_{M}(c\phi(\partial\sigma ))+1\leq 
diam_{M}(\phi(\partial\sigma ))+2+1\leq b+3$$
Thus there is a $\gamma$ such that $c\phi (\partial\sigma )\subset 
\gamma C_{0}$.
Let $g$ be the map corresponding to $\gamma^{-1}$ given by claim 2.

Now, for $x\in\partial \sigma $,
$$d_{K}(cg\phi(x),C_{0})\leq d_{K}(\gamma^{-1}c\phi(x),C_{0})+15 = 
d_{K}(c\phi(x),\gamma C_{0})+15=15$$
which means that $cg\phi(x)\in C_{1}$. Thus $g\phi(\partial \sigma 
)\subset c^{-1}(C_{1})$.
This implies that there is a map $\phi ':\sigma\rightarrow 
c^{-1}(E_{0})\cap \partial M$ extending $g\phi $. 
Extend now $\phi$ by defining $\phi =g^{-1}\phi '$.
Now, for $x\in\sigma$,
$$d_{K}(c\phi(x),\gamma E_{0})= d_{K}(cg^{-1}\phi'(x),\gamma E_{0})= 
d_{K}(\gamma^{-1}cg^{-1}\phi' (x),E_{0}) \leq d_{K}(cgg^{-1}\phi 
'(x),E_{0})+15 =15$$
consequently $c\phi(x)\in \gamma E_{1}$. Thus $\phi(\sigma )\subset 
c^{-1}(\gamma E_{1})$, and we get
(use {\bf (1)})
 $$diam_{M}(\phi (\sigma ))\leq diam_{M}(c^{-1}(\gamma E_{1}))\leq 
diam_{M}(\gamma E_{1})+2\leq diam_{K}(\gamma E_{1})+2= diam_{K}(E_{1})+2$$
and take $a=diam_{K}(E_{1})+2$, that depends only on $b$. This 
completes the proof of the claim.
\vspace{.2in}

{\bf Claim 4.} There is a constant $a$, such that the following holds. 
For any pair
 of finite simplicial complexes $L_{0}\subset L$, with $dim(L)\leq 
(n-1)$, 
 $L_{0}$ contains the vertices of $L$, and a $PL$
map $\psi : L_{0}\rightarrow \partial M$ with $diam_{M}(\psi 
(L_{0}\cap\sigma ))\leq 3$, $\sigma\in (L)^{1}$,  there is an extension $\psi:  
L\rightarrow \partial M$, with $diam_{M}(\psi (\sigma ))\leq a$, for all  
$\sigma\in L$.
\vspace{.2in}

Here $(L)^{1}$ is the one-skeleton of $L$. To prove claim 4, 
proceed as follows.
 First, extend $\psi$ to $L_{0}\cap (L)^{1}$. For this
just choose a geodesic between the endpoints of $\sigma^{1}\in 
(L)^{1}$. Note that
now $diam_{K}(\psi (\partial\sigma^{2}))\leq 6$, for every two-simplex 
$\sigma^{2}$.
 Apply claim 3 to $b=6$, and we get a constant $a_{6}$, and an 
extension of $\psi$ to
$L_{0}\cap (L)^{2}$, with $diam_{K}(\psi (\sigma^{2}))\leq a_{6}$. Then 
$diam_{K}(\psi (\partial\sigma^{3}))\leq 2a_{6}$, and apply claim 3 
again. We proceed
in the same way until we reach dimension $(n-1)$. This proves claim 4.
\vspace{.2in}

Now, let $N=a+2$, where $a$ is the constant from claim 4, and choose a 
base point $x_{0}\in \partial M $.

Recall that we assumed $(M,\partial M)$ to be  $PL$-homeomorphic 
to $({\Bbb{R}}^{n}_{+},{\Bbb{R}}^{n-1})$. Recall also that $(M,d_{M})$ is
a proper metric space.

Let $A$ denote the closed ball in $M$ with center $x_{0}$ and radius $N$.
Since $A$ is compact, there is a $PL$ $(n-1)$-ball ${\tilde D}\subset 
\partial M$, such that $A\cap\partial M\subset int\, \tilde D$.

Let ${\tilde {\rho}} :S\rightarrow \partial M $, where $S$ is a $(n-2)$-sphere, be 
a $PL$-embedding
such that ${\tilde{\rho}} (S)=\partial ({\tilde D})$. Because $x_{0}\in int\, {\tilde D}$
we have that ${\tilde{\rho}}\neq 0\in\pi_{n-2}(\partial M\setminus\{ x_{0}\})$

Remark also that  ${\tilde{\rho}}=0\in \pi_{n-2}(M\setminus A)$.
Thus there is a $PL$ extension $\rho$ of ${\tilde{\rho}}$ with
$\rho:  D\rightarrow  M\setminus A $,  where $D$ is a closed $(n-1)$-ball 
and $\partial D=S$.
Hence $$d_{M}(\rho (v), x_{0})\geq N,\,\,\,\, for\,\,\,\, all\,\,\,\, v\in D$$ 

and
$$\rho|_{S}\neq 0\in\pi_{n-2}(\partial M\setminus\{ x_{0}\}).$$
\vspace{.1in}

Subdivide $D$ so that $diam_{M}(\rho (\sigma))\leq 1$, for any simplex 
$\sigma \in D$. Denote by $D_{0}=\{ v_{1},...,v_{k}\}$ the zero-skeleton 
of $D$. Then  $v_{1},...,v_{k}$ are the {\it vertices} of $D$. Because of {\bf (2)}, for 
every $v_{i}$ we can select a $w_{i}\in \partial M$ such that 
$d_{M}(\rho(v_{i}),w_{i})\leq 1$  (if $\rho (v_{i})\in\partial M$
choose $w_{i}=\rho (v_{i})$). We have that if $w_{i},w_{j}$ correspond 
to the  vertices $v_{i}$,$v_{j}$ of a simplex $\sigma\in D$ then
$$d_{M}(w_{i},w_{j})\leq 
d_{M}(w_{i},\rho(v_{i}))+d_{M}(\rho(v_{i}),\rho(v_{j}))+d_{M}(\rho(v_{j}),w_{j})\leq 3,$$

\noindent and since $d_{M}(\rho (v), x_{0})\geq N$, for all $v\in D$, 
we get $$d_{M}(w_{i},x_{0})\geq d_{M}(\rho (v_{i}),x_{0})-d_{M}(\rho 
(v_{i}), w_{i})\geq N-1.$$
Define $\psi:S\cup D_{0}\rightarrow \partial M$, by $\psi (u)=\rho 
(u)$, for $u\in S$, and $\psi (v_{i})=w_{i}$, for $i=1,...k$.  Remark that 
$d_{M}(\psi (u),x_{0})\geq N-1$, for $u\in S \cup D_{0}$.

By claim 4 we can extend $\psi$ to a map 
$\psi: D\rightarrow \partial M$, with $diam_{M}(\psi (\sigma ))\leq a$, for 
all  $\sigma\in D$. Note that, by definition, $\psi |_{S}=\rho |_{S}$. Then, for $u\in D$, 
$$d_{M}(\psi(u),x_{0})\geq  d_{M}(\psi (v_{i}), x_{0}) - d_{M}(\psi (u),\psi (v_{i}))     \geq N-1-a
=(a+2)-1-a=1>0$$
where $v_{i}$ is a vertex in a simplex that contains $u$.

This implies that $\psi (D)\subset \partial M\setminus \{x_{0}\}$.
Hence $\psi |_{S}=0\in\pi_{n-2}(\partial M\setminus\{ x_{0}\} )$. But 
this is a contradiction because $\psi |_{S}=\rho |_{S}\neq 0\in\pi_{n-2}(\partial 
M\setminus \{x_{0}\})$.
This completes the proof of the proposition.

\end{document}